\documentclass[11pt, a4paper]{article}
\usepackage{amsmath,amssymb,latexsym,graphics,epsfig}
\usepackage{hyperref}
\usepackage{color}
\usepackage{amsthm}

\newtheorem{theo}{Theorem}[section]

\def\prova{{\boldmath  $Proof.$}\hskip 0.3truecm}

\def\final{\mbox{ \quad $\Box$}}

\def\A{\mbox{\boldmath $A$}}
\def\B{\mbox{\boldmath $B$}}

\def\Part{{\cal P}}

\def\Re{\mathbb R}

\def\dist{\mathop{\rm dist }\nolimits}

\def\excu{\mbox{$\varepsilon_u$}}

\def\e{\mbox{\boldmath $e$}}

\def\vecrho{\mbox{\boldmath $\rho$}}

\def\vecalpha{\mbox{\boldmath $\alpha$}}

\def\A{\mbox{\boldmath $A$}}

\def\E{\mbox{\boldmath $E$}}

\def\G{\Gamma}

\def\I{\mbox{\boldmath $I$}}

\def\J{\mbox{\boldmath $J$}}

\def\sp{\mbox{\rm sp}}

\def\ecc{\mathop{\rm ecc }\nolimits}

\def\sp{\mathop{\rm sp }\nolimits}

\begin{document}
\title{Pseudo-Distance-Regularised Graphs Are \\ Distance-Regular or
Distance-Biregular
}

\author{M.A. Fiol
\\ \\
{\small Universitat Polit\`ecnica de Catalunya}\\
{\small Departament de Matem\`atica Aplicada IV} \\
{\small Barcelona, Catalonia} \\
{\small e-mail:  {\tt fiol@ma4.upc.edu}}
}


\maketitle

\begin{abstract}
The concept of pseudo-distance-regularity around a vertex of a graph is a natural generalization, for non-regular graphs, of the standard distance-regularity around a vertex. In this note, we prove that a pseudo-distance-regular graph around each of its vertices is either distance-regular or distance-biregular. By using a combinatorial approach, the same conclusion  was reached by  Godsil and  Shawe-Taylor  for a distance-regular graph around each of its vertices. Thus, our proof, which is of an algebraic nature, can also be seen as an alternative demonstration of Godsil and Shawe-Taylor's theorem.
\end{abstract}


Distance-biregular graph.

\section{Introduction}
Distance-regularity around a vertex of a (regular) graph is the local analogue of distance-regularity. More precisely,
a graph $\G$ with vertex set $V$ is distance-regular around of a vertex $u$  if the distance partition of $V$ induced by $u$ is regular (see, for instance, Brouwer, Cohen and Neumaier \cite{bcn89}). In \cite{gst87},
Godsil and Shawe-Taylor  defined a distance-regularised graph as that being distance-regular around each of its vertices. The interest of these graphs relies on the fact that they are a common generalization of distance-regular graphs and generalized polygons.
The authors of \cite{gst87} used a combinatorial approach to prove that distance-regularised graphs are either distance-regular or distance-biregular. For some properties of the latter, see Delorme \cite{d94}.
More recently, Fiol, Garriga and Yebra \cite{fgy96} introduced the concept of pseudo-distance-regularity around a vertex, as a natural generalization for not necessarily regular graphs, of distance regularity around a vertex. In this note, we prove that the same conclusion obtained in \cite{gst87} can be reached when $\G$ is {\em pseudo-distance-regularized}; that is, pseudo-distance-regular graph around each of its vertices.
In fact, this was already obtained in \cite{fgy96}, but making strong use of the result in \cite{gst87}. Here we provide an independent direct proof, which is simple and of algebraic nature.
Thus, our contribution can be seen as an alternative demonstration of Godsil and Shawe-Taylor's theorem. Moreover, it turns out that the same conclusion of the theorem is obtained from the seemingly weaker condition of pseudo-distance-regularity.

\section{Preliminaries}
Let us first give some basic notation and results on which our proof is based. For more background on
graph spectra, distance-regular and pseudo-distance-regular graphs see, for instance, \cite{bi84,b93,bcn89,bh12,cds82,dkt12,f01,f02,fgy96}.

Throughout this note, $\G$ is a connected graph with vertex set $V=V(G)$, $n=|V|$ vertices, adjacency matrix $\A$ and spectrum
$
\sp \G = \sp \A = \{\lambda_0^{m_0},\lambda_1^{m_1},\dots,
\lambda_d^{m_d}\},
$
where the different eigenvalues of $\G$ are in decreasing order,
$\lambda_0>\lambda_1>\cdots >\lambda_d$, and the superscripts
stand for their multiplicities $m_i=m(\lambda_i)$, $i=0,1,\ldots,d$.
Then, as it is well known, $\lambda_0$, with multiplicity $m_0=1$,  coincides with the spectral radius of $\A$, and has a positive (column) eigenvector (the {\em Perron vector\/}) $\vecalpha$, which we normalize in such a way that $\|\vecalpha\|^2=n$. Let $\dist(u,v)$ represent the distance between the vertices $u,v\in V$. Then, $\G_i(u)=\{v|\dist(u,v)=i\}$, the {\em eccentricity} of a vertex $u$ is $\ecc(u)=\max\{\dist(u,v)|v\in V\}$, and the {\em diameter} of $\G$ is $D=\max\{\ecc(u)|u\in V\}$. For every $0\le i\le D$, the {\em distance-$i$ matrix} $\A_i$ has entries $(\A_i)_{uv}=1$ if $\dist(u,v)=i$, and $(\A_i)_{uv}=0$, otherwise. We also use the {\em weighted distance-$i$ matrix}, defined as $\A_i^*= \A_i\circ \J^*$, where $\J^*=\vecalpha \vecalpha^{\top}$ and `$\circ$' stands for the Hadamard product; that is, $(\A_i^*)_{uv}=\alpha_u\alpha_v$ if $\dist(u,v)=i$, and $(\A_i^*)_{uv}=0$, otherwise.

\subsection{Local spectrum and predistance polynomials}

Given a graph $\G$ with adjacency matrix $\A$ and spectrum as above, its idempotents
$\E_i$ correspond to the orthogonal projections onto the eigenspaces $\ker (\A-\lambda_i\I)$, $i=0,1,\ldots,d$. Their entries $m_{uv}(\lambda_i)=(\E_i)_{uv}$ are called the {\em (crossed) $uv$-local multiplicities} of $\lambda_i$. In particular, the diagonal elements $m_u(\lambda_i)=m_{uu}(\lambda_i)$ are the so-called {\em $u$-local multiplicities} of $\lambda_i$, because they satisfy properties similar to those of the (global) multiplicities $m(\lambda_i)$, but when $\G$ is ``seen" from the base vertex $u$. Indeed,
\begin{equation}\label{propert-locmul}
\textstyle a_{uu}^{(\ell)}=(\A^{\ell})_{uu}=\sum_{i=0}^d m_u(\lambda_i)\lambda_i^{\ell}, \quad\mbox{and}\quad \sum_{u\in V}m_u(\lambda_i)=m(\lambda_i).
\end{equation}
The {\em $u$-local spectrum} of $\G$ is constituted by the eigenvalues of $\A$, say $\mu_0,\mu_1,\ldots,\mu_{d_u}$, with nonzero $u$-local multiplicity. Then, it is known that the vector space spanned by the vectors $(\A^{\ell})_u$, $\ell\ge 0$, (that is, the $u$-th columns of the matrices $\A^{\ell}$, $\ell\ge 0$) has dimension $d_u$, and the eccentricity of $u$ satisfies $\ecc(u)\le d_u$. When $\ecc(u)= d_u$ we say that $u$ is {\em extremal} (for more details, see \cite{fgy96}).

An orthogonal base for such a vector space is the following.
The {\em $u$-local predistance polynomials} $p_0^u,p_1^u,\ldots, p_{d_u}^u$, $\deg p_i=i$, associated to a vertex $u$ of $\G$ with nonzero local multiplicities $m_u(\mu_i)$, $i=0,1,\ldots,d_u$, are a sequence of orthogonal polynomials with respect to the scalar product
$$
\langle f,g\rangle_u = (f(\A)g(\A))_{uu}= \sum_{i=0}^{d_u} m_u(\mu_i) f(\mu_i)g(\mu_i)= \sum_{i=0}^{d} m_u(\lambda_i) f(\lambda_i)g(\lambda_i),
$$
normalized in such a way that  $\|p_i^u\|_u^2=\alpha_u^2p_i^u(\lambda_0)$.
We notice that in \cite{fgy96} the normalization condition was $\|p_i^u\|_u^2=p_i^u(\lambda_0)$ but, although the theory remains unchanged, it seems more convenient to use the above.
Then, in particular, $p_0^u=\alpha_u^2$ and $p_1^u=\frac{\alpha_u^2\lambda_0}{\delta_u} x$. Indeed, they are orthogonal since $\langle 1,x \rangle_u= \sum_{i=0}^{d_u} m_u(\lambda_i) \lambda_i=0$, and the normalization condition is fulfilled:
\begin{itemize}
\item
$\|\alpha_u^2\|_u^2 = \alpha_u^4\sum_{i=0}^{d}  m_u(\lambda_i) =\alpha_u^4=\alpha_u^2p_0^u(\lambda_0)$.
\item
$\|\frac{\alpha_u^2\lambda_0}{\delta_u} x\|_u^2 = \frac{\alpha_u^4\lambda_0^2}{\delta_u^2}\sum_{i=0}^{d} m_u(\lambda_i) \lambda_i^2=  \frac{\alpha_u^4\lambda_0^2}{\delta_u}=\alpha_u^2 p_1^u(\lambda_0)$.
\end{itemize}

If $\G$ is $\delta$-regular, then $\alpha_u=1$, $\lambda_0=\delta$, and we have $p_0=1$ and $p_1=x$, which are the first distance polynomials for every distance-regular graph. More generally, and as expected, if $\G$ is distance-regular, the $u$-local predistance polynomials are independent of $u$ and become the {\em distance polynomials} $p_i$, $i=0,1,\ldots, D$, satisfying
\begin{equation}
\label{pi(A)=Ai}
p_i(\A)=\A_i,\qquad i=0,1,\ldots,D.
\end{equation}
As every sequence of orthogonal polynomials, the $u$-local predistance polynomials satisfy a three-term recurrence of the form
\begin{equation}\label{recur-pol}
 xp_i^u=b_{i-1}^*p_{i-1}^u+a_{i}^*p_{i}^u+c_{i+1}^*p_{i+1}^u, \qquad i=0,1,\ldots,d_u,
\end{equation}
where $b_{-1}^*=c_{d_u+1}^*=0$, and the other numbers $b_{i-1}^*$, $a_{i}^*$, and $c_{i+1}^*$ are the Fourier coefficients of $xp_i^u$ in terms of $p_{i-1}^u$, $p_{i}^u$, and $p_{i+1}^u$, respectively.


\subsection{Pseudo-distance-regularity around a vertex}
\label{sec: d-r conjunts}

Given a graph $\G$ as above,
consider the mapping $\vecrho: V \longrightarrow \Re^n$ defined by
$\vecrho(u)=\alpha_u\e_u$, where $\e_u$ is the coordinate vector.
Note that, since $\|\vecrho(u)\|=\alpha_u$, we can see $\vecrho$ as a function which assigns weights to the vertices of $\G$. In doing so we ``regularize" the graph, in
the sense that the {\it average weighted degree} of each vertex $u\in V$ becomes a
constant:
\begin{equation}
\label{regularize}
\textstyle \delta^*_u=\frac{1}{\alpha_u}\sum_{v\in\G(u)}\alpha_v=\lambda_0,
\end{equation}
where $\G(u)=\G_1(u)$.
Using these weights, we  consider the following concept. A partition $\Part$ of the vertex set
$V=V_1\cup\cdots\cup V_m$ is called {\it pseudo-regular} (or {\it pseudo-equitable})
whenever the {\it pseudo-intersection numbers}
\begin{equation}
\label{intersec-num}
b_{ij}^*(u)=\frac{1}{\alpha_u}\sum_{v\in\G(u)\cap V_j} \alpha_v,\qquad  i,j=0,1,\ldots,m
\end{equation}
do not depend on the chosen vertex $u\in V_i$, but only on the subsets $V_i$ and $V_j$.
In this case, such numbers are simply written as $b_{ij}^*$, and the $m\times m$ matrix
$\B^*=(b_{ij}^*)$ is referred to as  the {\it pseudo-quotient matrix} of $\A$ with respect to the
(pseudo-regular) partition $\Part$.
Pseudo-regular partitions were introduced by Fiol and Garriga \cite{fg99}, as a generalization of the so-called regular partitions, where the above
numbers are defined by $b_{ij}^*(u)=|\G(u)\cap V_j|$ for $u\in V_i$. A detailed study of regular partitions can be found in Godsil \cite{g93} and Godsil and  McKay \cite{gmk80}.

Let $u$ be a vertex of $\G$ with eccentricity $\ecc(u)=\excu$. Then $\G$ is said to be {\em pseudo-distance-regular around $u$}  (or {\em $u$-local pseudo-distance-regular\/}) if the {\em distance-partition} around $u$, that is
$V=C_0\cup C_1\cup\cdots\cup C_{\excu}$ where $C_i=\G_i(u)$ for $i=0,1,\ldots,\excu$, is pseudo-regular.
From the characteristics of the distance-partition, it is clear that its pseudo-quotient  matrix is tridiagonal $\B^*=(b_{ij}^*)$ with nonzero entries $c^*_i=b^*_{i-1,i}$, $a^*_i=b^*_{i,i}$ and $b^*_i=b^*_{i+1,i}$, $0\le i\le \excu$, with the convention $c^*_0=b^*_{\excu}=0$. By (\ref{regularize}), notice that
$a^*_i+b^*_i+c^*_i=\lambda_0$ for $i=0,1\ldots,\excu$. These parameters are called the {\em $u$-local $($pseudo-$)$intersection numbers}.
In \cite{fgy96}, it was shown that local pseudo-distance regularity is a generalization of distance-regularity around a vertex. Indeed, if $\G$ is distance-regular around $u\in V$, with intersection numbers $a_i,b_i,c_i$, then the entries of the Perron vector $\vecalpha$ have a constant value, say $\alpha^i$,  on each of the sets $\G_i(u)$, $i=0,1,\ldots,\excu$, and $\G$ turns out to be pseudo-distance-regular around $u$, with $u$-local intersection numbers
\begin{equation}
\label{abc*-abc}
a_i^*=a_i,\quad b_i^*=\frac{\alpha^{i+1}}{\alpha^i}b_i,\quad c_i^*=\frac{\alpha^{i-1}}{\alpha^i}c_i, \qquad i=0,1,\ldots,\excu.
\end{equation}
Conversely, when the eigenvector $\vecalpha$ of a pseudo-distance-regular graph $\G$ exhibits such a regularity (which is the case for all $u$ when $\G$ is regular or bipartite biregular), we have that $\G$ is also distance-regular around $u$ with intersection parameters given again by (\ref{abc*-abc}).

As happens for distance-regular graphs, the existence of the so-called {\em $u$-local distance polynomials}, satisfying the ``local version" of (\ref{pi(A)=Ai}), guarantees that $\G$ is pseudo-distance-regular around $u$.

\begin{theo}[\cite{fgy96}]
\label{teo-charac-pdrg-pol}
Let $\G$ be a graph having a vertex $u$
 with eccentricity $\excu$. Then, $\G$ is pseudo-distance-regular around $u$ if and only if the $u$-local predistance polynomials satisfy
$$
(p_i^u(\A))_u=(\A_i^*)_u,\qquad i=0,1,\ldots,\excu.
$$
Moreover, if this is the case, $u$ is extremal, $\excu=d_u$, and the $u$-local intersection numbers $a_i^*$ ,$b_i^*$ and $c_i^*$ coincide with the Fourier coefficients of the recurrence (\ref{recur-pol}).
\end{theo}
By this result, note that the $u$-local intersection numbers are univocally determined, through the $u$-local predistance polynomials, by the $u$-local spectrum.

\section{The proof}

Now we are ready to give the algebraic proof that a graph which is pseudo-distance-regular around each of its vertices is either distance-regular or distance-biregular.

\begin{theo}
\label{main-theo}
Every  pseudo-distance-regularized graph $\G$ is either distance-regular or distance-biregular.
\end{theo}
\prova
Let $v,w$ be two vertices adjacent to a vertex $u$. Then, by Theorem \ref{teo-charac-pdrg-pol},
$$
\alpha_u\alpha_v  =  (p_1^u(\A))_{uv}=\frac{\alpha_u^2\lambda_0}{\delta_u}(\A)_{uv}
=\frac{\alpha_u^2\lambda_0}{\delta_u}(\A)_{uw}=
(p_1^u(\A))_{uw}
=\alpha_u\alpha_w,
$$
and we infer that $\alpha_v=\alpha_w$ since $\alpha_u>0$. Hence, since $\G$ is connected, all vertices at even (respectively odd) distance from $u$ have component $\alpha_u$ (respectively $\alpha_v$), and $\G$ is either bipartite biregular if $\alpha_u\neq \alpha_v$, or regular otherwise.

Moreover, by using the orthogonal decomposition of $x^{\ell}$ in terms of the base $\{p_i^u\}_{0\le i\le d_u}$, the number of $\ell$-walks between two adjacent vertices $u,v$ is
\begin{equation}\label{l-walks-u}
(\A^{\ell})_{uv}=\sum_{i=0}^{d_u} \frac{\langle x^{\ell},p_i^u\rangle_u}{\|p_i^u\|_u^2}(p_i^u(\A))_{uv}
= \frac{\langle x^{\ell},p_1^u\rangle_u}{\alpha_u^2 p_1^u(\lambda_0)}\alpha_u\alpha_v
= \frac{\alpha_v}{\alpha_u}\frac{1}{\lambda_0}\sum_{i=0}^d m_u(\lambda_i)\lambda_i^{\ell+1}
\end{equation}
and, similarly,
\begin{equation}\label{l-walks-v}
(\A^{\ell})_{uv}=
\frac{\alpha_u}{\alpha_v}\frac{1}{\lambda_0}\sum_{i=0}^d m_v(\lambda_i)\lambda_i^{\ell+1}
\end{equation}
Hence, from (\ref{propert-locmul}), (\ref{l-walks-u}), and (\ref{l-walks-v}) we have that, given the local multiplicities of $u$, those of $v$ are uniquely determined from the system
\begin{eqnarray*}
\textstyle \sum_{i=0}^d m_v(\lambda_i) &=& 1 \\
\textstyle \sum_{i=0}^d m_v(\lambda_i)\lambda_i^{r} &=& \textstyle \frac{\alpha_u^2}{\alpha_v^2}\sum_{i=0}^d m_u(\lambda_i)\lambda_i^{r}, \qquad r=1,2,\ldots, d,
\end{eqnarray*}
with $d$ equations, $d$ unknowns
$m_v(\lambda_i)$ ($i=1,2,\ldots,d$), and Vandermonde determinant.
Then, we have only two possible cases:
\begin{enumerate}
\item
If $\G$ is regular, then $\alpha_u=\alpha_v=1$ for every $u,v\in V$ and, hence, every vertex has the same local spectrum (this is known to be equivalent to say that $\G$ is {\em walk-regular} \cite{gmk80}), and the same $u$-local distance-polynomials. Thus, $\G$ is distance-regular.
\item
If $\G$ is bipartite $(\delta_1,\delta_2)$-biregular, then $\alpha_u=\sqrt{(\delta_1+\delta_2)/2\delta_2}$ for every $u\in V_1$ and $\alpha_v=\sqrt{(\delta_1+\delta_2)/2\delta_1}$ for every $v\in V_2$. In this case, every vertex of $V_1$ has the same local spectrum and the same holds for every vertex in $V_2$ (in this case, we could say that $\G$ is {\em walk-biregular}). Thus, the sequence of local $u$-predistance polynomials only depend on the partite set where $u$ belongs to and, consequently, $\G$ is distance-biregular.
\end{enumerate}
This completes the proof.
\final

As mentioned above, notice that, since every regular or biregular graph which is pseudo-distance-regular around a vertex is also distance-regular around such a vertex, Theorem \ref{main-theo} implies the result of Godsil and Shawe-Taylor in \cite{gst87}.

\noindent {\bf Acknowledgments.}
Research supported by the Ministerio de Educaci\'on y
Ciencia (Spain) and the European Regional Development Fund under
project MTM2011-28800-C02-01, and by the Catalan Research Council
under project 2009SGR1387.

\vskip -1cm



\begin{thebibliography}{99}


\bibitem{bi84}
E. Bannai and T. Ito, \emph{Algebraic Combinatorics I: Association Schemes}, Benjamin/Cummings,
London, 1974, 1984.

\bibitem{b93}
N. Biggs, \emph{Algebraic Graph Theory}, Cambridge University
Press, Cambridge, 1974, second edition, 1993.

\bibitem{bcn89}
A.E. Brouwer, A.M. Cohen, and A. Neumaier, \emph{Distance-Regular Graphs},
Springer-Verlag, Berlin-New York, 1989.

\bibitem{bh12}
A.E. Brouwer and W.H. Haemers, \emph{Spectra of graphs},
Springer, 2012; available online at
\url{http://homepages.cwi.nl/~aeb/math/ipm/}.


\bibitem{cds82}
D. M. Cvetkovi\'c, M. Doob and H. Sachs, Spectra of Graphs. Theory and Application,
VEB Deutscher Verlag der Wissenschaften, Berlin, second edition, 1982.


\bibitem{dkt12}
E.R. van Dam, J.H. Koolen, and H. Tanaka, Distance-regular graphs,
manuscript (2012), available online at
\url{http://lyrawww.uvt.nl/~evandam/files/drg.pdf}.

\bibitem{d94}
C. Delorme, Distance biregular bipartite graphs,
{\em Europ. J. Combin.} {\bf 15} (1994), 223--238.

\bibitem{f01}
M.A. Fiol, On pseudo-distance-regularity.
\emph{Linear Algebra Appl.} {\bf 323} (2001), 145--165.

\bibitem{f02}
M.A. Fiol, Algebraic characterizations of distance-regular graphs,
\emph{Discrete Math.} {\bf 246} (2002), 111--129.



\bibitem{fg99}
M.A. Fiol and E. Garriga, On the algebraic theory
of pseudo-distance-regularity around a set, {\it Linear
Algebra Appl.} {\bf 298} (1999), 115--141.



\bibitem{fgy96}
M.A. Fiol, E. Garriga, and J.L.A. Yebra,
Locally pseudo-distance-regular graphs,
{\em J.  Combin. Theory Ser. B} {\bf 68} (1996), 179--205.

\bibitem{g93}
C.D. Godsil, {\it Algebraic Combinatorics}, Chapman and Hall, NewYork, 1993.

\bibitem{gmk80}
C.D. Godsil and B.D. McKay,
Feasibility conditions for the existence of walk-regular graphs,
{\em Linear Algebra Appl.} {\bf 30} (1980) 51--61.

\bibitem{gst87}
C.D. Godsil and J. Shawe-Taylor,
Distance-regularised graphs are distance-regular or distance-biregular,
{\em J. Combin. Theory Ser. B} {\bf 43} (1987) 14--24.






\end{thebibliography}
\end{document}